\documentclass[12pt]{amsart}
\usepackage{amsfonts}
\usepackage{amssymb}
\usepackage{amsmath}
\usepackage{amsthm}
\usepackage{cases}
\usepackage{graphicx}
\usepackage{fullpage}

\theoremstyle{plain}
\newtheorem{thm}{Theorem}[section]
\newtheorem{lem}[thm]{Lemma}
\theoremstyle{definition}
\newtheorem{defn}[thm]{Definition}
\newtheorem{ex}[thm]{Example}

\numberwithin{equation}{section}

\makeatletter

\begin{document}
\title[A nonlocal fractional boundary value problem]{Positive solutions of a nonlocal Caputo fractional BVP}\thanks{Published in: Dynamic Systems and Applications 23 (2014) 715--722.}%
\date{}


\subjclass[2010]{Primary 34A08, secondary 34B10, 34B18.}%
\keywords{Positive solution, nonlocal boundary conditions, fractional equation, fixed point index, cone.}%

\author[A. Cabada]{Alberto Cabada$^*$}\thanks{$^*$Partially supported by FEDER and Ministerio de Educaci\'on y Ciencia, Spain, project MTM2010-15314.}
\address{Alberto Cabada, Departamento de An\'alise Ma\-te\-m\'a\-ti\-ca, Facultade de Matem\'aticas, 
Universidade de Santiago de Com\-pos\-te\-la, 15782 Santiago de Compostela, Spain}%
\email{alberto.cabada@usc.es}%

\author[G. Infante]{Gennaro Infante}
\address{Gennaro Infante, Dipartimento di Matematica ed Informatica, Universit\`{a} della
Calabria, 87036 Arcavacata di Rende, Cosenza, Italy}%
\email{gennaro.infante@unical.it}%

\begin{abstract}
We discuss the existence of multiple positive solutions for a nonlocal fractional problem recently considered by Nieto and Pimentel. Our approach relies on classical fixed point index.
\end{abstract}

\maketitle

\centerline{\it
Dedicated to Professor John R. Graef on the occasion of his seventieth birthday.
}
\section{Introduction}
Recently, Nieto and Pimentel~\cite{nieto-pim} studied the existence of positive solutions of the nonlocal fractional boundary value problem (BVP)
 \begin{gather}\label{niepim}
 \begin{aligned}
 {}^C\!D^{\alpha} u(t)+f(t, u(t))=0, \ t\in(0, 1), \\
 u'(0)=0,\ \beta{}^C\!D^{\alpha-1}u(1)+u(\eta)=0,
  \end{aligned}
 \end{gather}
where $1<\alpha\leq 2$, ${}^C\!D^{\alpha}$ denotes the Caputo fractional derivative of order $\alpha$, $\beta>0$, $0\leq \eta\leq 1$ and $f$ is continuous. The reason, given in~\cite{nieto-pim}, for studying the BVP~\eqref{niepim} is that it is seen as a mathematical generalisation of the BVP
\begin{gather}\label{therm}
 \begin{aligned}
u''(t)+f(t, u(t))=0, \ t\in(0, 1), \\
 u'(0)=0,\ \beta u'(1)+u(\eta)=0,
  \end{aligned}
   \end{gather}
that was studied by Infante and Webb~\cite{gijwnodea}, who were motivated by the previous work of Guidotti and Merino~\cite{guimer}.
The BVP~\eqref{therm} can be used as a model for 
heated bar of length~1 with a thermostat,
where a controller at $t=1$ adds or removes heat according to the
temperature detected by a sensor at $t=\eta$. Heat-flow problems of this type have been studied recently, see for example \cite{Fan-Ma, gi-poit, gi-caa, gijwems, Kar-Pal, pp-gi-pp-aml-06, jwpomona, jwwcna04, jw-narwa} and references therein.

Here we discuss the existence of multiple positive solutions for the nonlocal BVP 
 \begin{gather}\label{fbvp}
 \begin{aligned}
 {}^C\!D^{\alpha} u(t)+f(t, u(t))=0, \ t\in(0, 1), \\
 u'(0)+\lambda[u]=0,\ \beta{}^C\!D^{\alpha-1}u(1)+u(\eta)=0,
  \end{aligned}
 \end{gather}
where
$\lambda[\cdot]$ is a functional given by
$$
\lambda[u]=\Lambda_0+\int_0^1 u(s)\,d\Lambda(s),
$$
involving a Stieltjes integral.
This type of BCs includes as special cases
$$
\lambda[u]=\sum_{i=1}^{m} \lambda_{i}u(\xi_{i})\quad (m\text{-point problems})
$$ 
and 
$$
\lambda[u]=\int_{0}^{1} {\lambda}(s)u(s)\,ds\quad (\text{continuously distributed cases}).
$$
Multi-point and integral BCs are widely studied objects, see, for
example, Karakostas and Tsamatos \cite{kttmna,ktejde}, Ma \cite{ma}, 
Ntouyas \cite{sotiris}, 
Webb \cite{jw-seville,jw-poitiers}, Henderson and co-authors
\cite{hnp}, 
Infante and Webb \cite{gijwems,jwginodeaII},  Zima \cite{miraJMAA04}. 

We mention that John Graef and co-authors have actively contributed to the study of nonlocal and fractional problems with interesting papers, for some of their recent works see~\cite{Graef1, Graef2, Graef3, Graef4}.

In this note, we use the methodology developed in \cite{gijwems}, that is to rewrite the BVP \eqref{fbvp} as a
\emph{perturbed} Hammerstein integral equation of the form
\begin{equation}\label{perhaminteq}
u(t)=\gamma(t){\lambda}[u]+\int_{0}^{1}k(t,s)f(s,u(s))\,ds,
\end{equation}
and use the classical \emph{theory of fixed point index}, for example see \cite{amann, guolak}, in order to
gain the existence of positive solutions of \eqref{perhaminteq}, by
working in a suitable cone of positive functions. The existence of positive solutions of \eqref{perhaminteq} provide the existence of \emph{positive solutions} of
the BVP~\eqref{fbvp}.

\section{Preliminaries}
We firstly recall the definition of the Caputo derivative. For its properties we refer to the books~\cite{anast, dieth, pod, samk}. 

\begin{defn}
For a function $y:[0, +\infty)\to \mathbb{R}$,
the  Caputo  derivative of fractional order
  $\alpha>0$  is given by
$$
{}^C\!D^\alpha
y(t)=\frac{1}{\Gamma(n-\alpha)}\int^t_0\frac{y^{(n)}(s)}{(t-s)^{\alpha+1-n}}\,ds, \quad n=[\alpha] + 1,
$$
where $\Gamma$ denotes the Gamma function, that is
$$\Gamma (s)=\int_{0}^{+\infty}x^{s-1}e^{-x}dx,$$
and $[\alpha]$ denotes the integer part of a number $\alpha$.
\end{defn}

We now recall some results from \cite{gijwems}, regarding the existence of multiple positive solutions of perturbed Hammerstein integral equations \begin{equation}\label{T}
u(t)=\gamma(t){\lambda}[u]+\int_{0}^{1}k(t,s)f(s,u(s))\,ds:=Tu(t),
\end{equation}

We point out that a cone $K$ in a Banach space $X$ is a closed,
convex set such that $\lambda x\in K$ for all $x \in K$ and $\lambda\geq
0$, and $K\cap (-K)=\{0\}$.

We work in the space of continuous functions $C[0,1]$ endowed with the usual supremum norm and we look for fixed points of $T$ in the following cone of non-negative functions
\begin{equation}\label{eqcone}
K=\{u\in C[0,1], u\geq 0: \min_{t \in [a,b]}u(t)\geq c \|u\|\},
\end{equation}
with $c$ a positive constant related to the kernel $k$ and the function $\gamma$.

The cone \eqref{eqcone} was first used by Krasnosel'ski\u\i{}, see e.g. \cite{krzab}, and D.~Guo, see e.g. \cite{guolak}, and then used by many authors.
 
The following assumptions on the terms that occur in \eqref{T} are a special case of the ones in~\cite{gijwems}.

\begin{itemize}\item  $f:[0,1]\times[0,\infty)\to [0,\infty)$ is continuous.
\item  $k:[0,1]\times [0,1]\to [0,\infty)$ is continuous.{}
\item There exist a function  $\Phi:[0,1]\to
[0,\infty)$,  $\Phi \in  L^1[0,1]$, an interval $[a,b]\subset [0,1]$ and a constant $c_{1} \in (0,1]$ such that
\begin{align*}
  k(t,s)\leq \Phi(s) \text{ for }& t,s \in [0,1] \text{ and}\\
  c_{1} \Phi(s)\leq k(t,s) \text{ for }& t \in [a,b] \text{ and }s \in
  [0,1].
\end{align*}{}
\item  ${\lambda}$ is an affine
functional given by
\begin{equation}
\label{e-lambda}
\lambda[u]= \Lambda_{0}+\int_{0}^{1} u(s)\,d\Lambda(s),
\end{equation}
where $ \Lambda_{0}\geq 0$ and $d\Lambda$ is a \emph{positive} Stieltjes measure with $\Lambda_1:=\int_{0}^{1}d\Lambda(s)<\infty.$
\item $\gamma :[0,1]
\to [0,\infty)$ is continuous, there exists a constant $c_{2}\in (0,1]$
such that 
$$
\gamma (t) \geq c_{2}\|\gamma\| \;\text{ for } t\in [a,b].
$$
and
$$
\tilde{\lambda}[\gamma]:=\int_{0}^{1} \gamma (t) \,d\Lambda(t)<1,
$${}
\end{itemize}
The assumptions above enable us to use the cone \eqref{eqcone} with $c=\min\{c_1,c_2\}$. A routine argument shows that $T$ maps $K$ into $K$ and is compact.

We make use, for our index calculations, of the following open bounded sets (relative to $K$):
$$K_{\rho}=\{u\in K: \|u\|<\rho\},\quad V_{\rho}=\{u \in K: \displaystyle{\min_{t\in[a,b]}}u(t)<\rho\}.$$
These sets have the key property that 
$$K_{\rho}\subset V_{\rho}\subset K_{\rho/c}.$$

The first Lemma is a special case of Lemma 2.4 of \cite{gijwems} and ensures that, for a suitable $\rho>0$, the fixed point index is 0 on the set $V_{\rho}$.
\begin{lem}
Assume that
\begin{enumerate}
\item[$(\mathrm{I}_{\protect\rho }^{0})$]
 there exists $\rho> 0$ such that for some ${\lambda}_{0}\geq 0$
\begin{equation}
\label{e-i=0}
{\lambda}[u]\geq {\lambda}_{0} \rho \;\text{ for } u\in \partial
V_{\rho},
\quad
c_{2}\|\gamma\| \lambda_{0} +f_{{\rho},{\rho / c}}
\cdot\frac{1}{M}> 1,
\end{equation}
where
\begin{align*}
 f_{{\rho},{\rho / c}}=\inf\Bigl\{\frac{f(t,u)}{\rho} :(t,u) \in [a,b]\times  [\rho,\rho / c]\Bigr\}\ \text{and}\
\frac{1}{M}=\inf_{t\in [a,b]}\int_{a}^{b} k(t,s)  \,ds.
\end{align*}
\end{enumerate}
Then $i_{K}(T,V_{\rho})=0$.

\end{lem}

The next Lemma is a special case of Lemma 2.6 of \cite{gijwems} provides a condition that yields, for a suitable $\rho>0$, that the index is 1 on a set $K_{\rho}$.

\begin{lem}
\label{ind1} 
Assume that 
\begin{enumerate}
\item[$(\mathrm{I}_{\protect\rho }^{1})$] 
there exists $\rho> 0$ such that  
\begin{equation}\label{EqH}
  \frac{\Lambda_{0}\| \gamma \|}{\rho(1-\tilde{\lambda}[\gamma])}
  +\Bigl( \frac{\| \gamma \|}{1-\tilde{\lambda}[\gamma]}\int_{0}^{1}\mathcal{K}(s) \,ds
  +\frac{1}{m} \Bigl) f^{0,\rho} < 1,
\end{equation}
where
\begin{align*}
\mathcal{K}(s)&=\int_{0}^{1} k(t,s) \,d\Lambda(t),\\
f^{0,{\rho}}&=\sup \Bigl\{\frac{f(t,u)}{\rho}:(t,u) \in [0,1]\times  [0,\rho]\Bigr\}\ \text{and}\\ \frac{1}{m}&=\sup_{t\in [0,1]}
\int_{0}^{1} k(t,s)  \,ds .
\end{align*}
\end{enumerate}
Then $i_{K}(T,K_{\rho})=1$.
\end{lem}
The two Lemmas above give the following result on the existence of multiple positive solutions for
Eq.~\eqref{T}. The proof follows from the properties of fixed point index and is omitted.

\begin{thm}
The integral equation \eqref{T} has at least one non-zero solution
in $K$ if any of the following conditions hold.

\begin{enumerate}

\item[$(S_{1})$] There exist $\rho _{1},\rho _{2}\in (0,\infty )$ with $\rho
_{1}/c<\rho _{2}$ such that $(\mathrm{I}_{\rho _{1}}^{0})$ and $(\mathrm{I}_{\rho _{2}}^{1})$ hold.

\item[$(S_{2})$] There exist $\rho _{1},\rho _{2}\in (0,\infty )$ with $\rho
_{1}<\rho _{2}$ such that $(\mathrm{I}_{\rho _{1}}^{1})$ and $(\mathrm{I}%
_{\rho _{2}}^{0})$ hold.
\end{enumerate}
The integral equation \eqref{T} has at least two non-zero solutions in $K$ if one of
the following conditions hold.

\begin{enumerate}

\item[$(S_{3})$] There exist $\rho _{1},\rho _{2},\rho _{3}\in (0,\infty )$
with $\rho _{1}/c<\rho _{2}<\rho _{3}$ such that $(\mathrm{I}_{\rho
_{1}}^{0}),$ $(
\mathrm{I}_{\rho _{2}}^{1})$ $\text{and}\;\;(\mathrm{I}_{\rho _{3}}^{0})$
hold.

\item[$(S_{4})$] There exist $\rho _{1},\rho _{2},\rho _{3}\in (0,\infty )$
with $\rho _{1}<\rho _{2}$ and $\rho _{2}/c<\rho _{3}$ such that $(\mathrm{I}%
_{\rho _{1}}^{1}),\;\;(\mathrm{I}_{\rho _{2}}^{0})$ $\text{and}\;\;(\mathrm{I%
}_{\rho _{3}}^{1})$ hold.
\end{enumerate}
The integral equation \eqref{T} has at least three non-zero solutions in $K$ if one
of the following conditions hold.

\begin{enumerate}
\item[$(S_{5})$] There exist $\rho _{1},\rho _{2},\rho _{3},\rho _{4}\in
(0,\infty )$ with $\rho _{1}/c<\rho _{2}<\rho _{3}$ and $\rho _{3}/c<\rho
_{4}$ such that $(\mathrm{I}_{\rho _{1}}^{0}),$ $(\mathrm{I}_{\rho _{2}}^{1}),\;\;(\mathrm{I}%
_{\rho _{3}}^{0})\;\;\text{and}\;\;(\mathrm{I}_{\rho _{4}}^{1})$ hold.

\item[$(S_{6})$] There exist $\rho _{1},\rho _{2},\rho _{3},\rho _{4}\in
(0,\infty )$ with $\rho _{1}<\rho _{2}$ and $\rho _{2}/c<\rho _{3}<\rho _{4}$
such that $(\mathrm{I}_{\rho _{1}}^{1}),\;\;(\mathrm{I}_{\rho
_{2}}^{0}),\;\;(\mathrm{I}_{\rho _{3}}^{1})$ $\text{and}\;\;(\mathrm{I}%
_{\rho _{4}}^{0})$ hold.
\end{enumerate}
\end{thm}

\section{The nonlocal BVP}
A direct calculation shows that solution of the linear equation $${}^C\!D^{\alpha} u(t)+y(t)=0,$$
 under the BCs 
$$
u'(0)+\lambda[u]=0,\ \beta{}^C\!D^{\alpha-1}u(1)+u(\eta)=0,
$$
can be written in the form
\begin{align*}
u(t)=&\Bigl(  \frac{\beta}{\Gamma(3-\alpha)} + \eta - t \Bigr){\lambda}[u] +\beta \int_{0}^{1}y(s)ds\\
&+\int_{0}^{\eta}\frac{(\eta-s)^{\alpha-1}}{\Gamma(\alpha)}y(s)ds -\int_{0}^{t}
\frac{(t-s)^{\alpha-1}}{\Gamma(\alpha)}y(s)ds.
\end{align*}
Therefore the solution of the BVP \eqref{fbvp} is
$$
u(t)=\gamma(t){\lambda}[u]
+\int_{0}^{1}k(t,s)f(s,u(s))ds
$$
where 
$$ \gamma (t)= \Bigl(  \dfrac{\beta}{\Gamma(3-\alpha)} + \eta - t \Bigr)$$ 
and
$$
k(t,s)=\beta +\frac{1}{\Gamma(\alpha)}\begin{cases} (\eta-s)^{\alpha-1},\ &s\leq \eta\\0,\ &s>\eta
\end{cases}
-\frac{1}{\Gamma(\alpha)}\begin{cases}(t-s)^{\alpha-1},\ &s\leq t\\ 0,\ &s>t.
\end{cases}
$$
Here we focus on the case $$\beta\Gamma(\alpha) >(1-\eta)^{\alpha-1}, \ \beta>(1-\eta)\Gamma(3-\alpha),$$
where $[a,b]$ can be chosen equal to $[0,1]$.

Upper and lower bounds for $k(t,s)$ were given in \cite{nieto-pim} as follows:
\begin{equation*}
\Phi(s)=\frac{\beta\Gamma(\alpha)+\eta^{\alpha-1}}{\Gamma(\alpha)},\quad  c_1=\frac{\beta\Gamma(\alpha)-(1-\eta)^{\alpha-1}}{\beta\Gamma(\alpha)+\eta^{\alpha-1}}.
\end{equation*}
Furthermore, by direct computation, we have
\begin{equation*}
\| \gamma \|=\eta+\frac{\beta}{\Gamma(3-\alpha)},\quad  c_2=\frac{\beta+(\eta-1)\Gamma(3-\alpha)}{\beta+\eta\Gamma(3-\alpha)}.
\end{equation*}

Hence we work with the cone
$$
K=\{u\in C[0,1]: \min_{ t \in [0,1]} u(t)\geq c \|u\|\},
$$
where
\begin{equation}\label{c-bcB}
c=\min\Bigl\{\frac{\beta\Gamma(\alpha)-(1-\eta)^{\alpha-1}}{\beta\Gamma(\alpha)+\eta^{\alpha-1}}, \frac{\beta+(\eta-1)\Gamma(3-\alpha)}{\beta+\eta\Gamma(3-\alpha)}\Bigr\}.
\end{equation}
\begin{ex}
Consider the BVP
 \begin{gather}
 \begin{aligned}
 {}^C\!D^{\alpha} u(t)+f(t,u(t))=0, \ t\in(0, 1),& \\
 u'(0)+\lambda u(\xi)=0,\ \beta{}^C\!D^{\alpha-1}u(1)+u(\eta)=0,& \ {\xi,\eta}\in [0,1].
  \end{aligned}
 \end{gather}

For this BVP  we may take in \eqref{e-lambda} $\Lambda_{0}=0$ and  $d\Lambda(s)$ the Dirac measure of
weight $\lambda >0$ at $\xi$.
A direct calculation gives
\begin{equation*}
m=\frac{\Gamma(\alpha+1)}{\beta\Gamma(\alpha+1)+\eta^\alpha}\quad \text{and}\quad \frac{1}{M}
=\frac{\Gamma(\alpha+1)}{\beta\Gamma(\alpha+1)+\eta^\alpha-1}.
\end{equation*}

Here we may take in \eqref{e-i=0} $\lambda_{0}=\lambda$ since, for $u\in \partial
V_\rho$, we have
$$
{\lambda}[u]= \lambda u(\xi)\geq \lambda \rho.
$$
For these BCs, \eqref{e-i=0} reads
\begin{equation}\label{gind03}
 \frac{\lambda(\beta+(\eta-1)\Gamma(3-\alpha))}{\Gamma(3-\alpha)}+f_{{\rho},{\rho / c}}
\cdot\frac{1}{M}>1,
\end{equation}
and so, $i_{K}(T,V_{\rho})=0$.

From Lemma \ref{ind1}, we have that $i_{K}(T,K_{\rho})=1$ if $\tilde{\lambda}[\gamma] <1$ and
\begin{equation}\label{EqH2}
\Bigl( \frac{\beta+\eta\Gamma(3-\alpha)}{(1-\tilde{\lambda}[\gamma])\Gamma(3-\alpha)}\int_{0}^{1}\mathcal{K}(s)\,ds
  +\frac{1}{m} \Bigl) f^{0,\rho} < 1.
\end{equation}
So we need
$$
\tilde{\lambda}[\gamma]=\int^{1}_{0} \gamma (t) \,d\Lambda(t)=\lambda \gamma
(\xi)=\lambda\Bigl(  \dfrac{\beta}{\Gamma(3-\alpha)} + \eta - \xi \Bigr)<1.
$$
Since $\mathcal{K}(s)=\lambda k(\xi,s)$, we obtain 
$$
\int_{0}^{1}\mathcal{K}(s)\,ds=\lambda\int_{0}^{1} k(\xi,s)\,ds=
\lambda \Bigl( \beta+\frac{\eta^\alpha}{\Gamma(\alpha+1)}-\frac{\xi^\alpha}{\Gamma(\alpha+1)} \Bigr).
$$
Note that all the numbers in \eqref{gind03} and \eqref{EqH2} can be computed. For example the choice
of $$\alpha= 3/2,\ \beta = 4/5,\ \eta = 3/4,\ \xi= 1/4,\ \lambda= 1/2$$
gives $c=0.132$. Then the $i_{K}(T,V_{\rho})=0$ condition needs
$$
f_{{\rho},{\rho / c}}> 0.218
$$
and
$i_{K}(T,K_{\rho})=1$
requires $$f^{0,\rho} < 1.255.$$
\end{ex}
\section*{Acknowledgments}
The authors would like to thank K. Lan for pointing out the recent paper \cite{lanfrac}, which contains interesting remarks on the relation between the solutions of a fractional BVP and the solutions of the corresponding integral equation. 

\end{document}